\documentclass[12pt]{amsart}
\usepackage{mathptmx}
\usepackage{amsmath}
\usepackage{amssymb}
\usepackage{array}
\usepackage{geometry}
\usepackage[bookmarks=true,colorlinks=true, pdfstartview=FitV, linkcolor=black, citecolor=blue, urlcolor=black]{hyperref}
\usepackage{color}
\definecolor{DarkRed}{rgb}{0.55,.00,0.2}
\definecolor{DarkGrey}{rgb}{0.35,.35,0.35}

\theoremstyle{definition}

\theoremstyle{remark}

\numberwithin{equation}{section}

\begin{document}

\title{On the Laguerre   fractional integro-differentiation}

\author{S. Yakubovich}
\address{Department of Mathematics, Fac. Sciences of University of Porto,Rua do Campo Alegre,  687; 4169-007 Porto (Portugal)}
\email{ syakubov@fc.up.pt}



\maketitle

\markboth{\rm \centerline{ S.  YAKUBOVICH}}{}
\markright{\rm \centerline{Laguerre fractional integro-differentiation}}

\begin{abstract} A fractional power interpretation of the Laguerre derivative $(DxD)^\alpha,\ D\equiv {d\over dx} $ is discussed.  The corresponding  fractional integrals are   introduced. Mapping and semigroup properties, integral representations and Mellin transform analysis are presented.   A relationship with the Riemann-Liouville  fractional integrals is demonstrated.   Finally,  a second kind integral equation of the Volterra-type, involving the Laguerre fractional integral is  solved in terms of the double hypergeometric type series as the resolvent kernel. 

\bigskip
\end{abstract}

\section{Introduction and auxiliary results}

The main goal of the present paper is to propose  constructions of the fractional integro-differentiation, which represent  fractional powers of the so-called Laguerre derivative $\theta\equiv DxD$ \cite{Dattoli},  where $D\equiv {d\over dx} $ is the differential operator.   The crucial feature of the Laguerre derivative is that its  integer powers satisfy   the following 
Viskov-type identity  (see \cite{Viskov})

$$\theta^n= \left( DxD\right)^n = D^n x^n D^n,\quad  n \in\mathbb{N}_{0}.\eqno(1.1)$$
Therefore it has a reason to give an interpretation of the fractional power $\left( DxD\right)^{\alpha}$ for arbitrary positive  values of $\alpha$.   

It is known \cite{Samka}, Section 5.1 that the classical left- and right-sided  Riemann-Liouville fractional integrals of order $\alpha >0$ on the half-axis ${\mathbb R}_+$ are defined, respectively,

$$\left(I_{0+}^\alpha f\right)(x) = {1\over \Gamma(\alpha)} \int_0^x (x-t)^{\alpha-1} f(t) dt,\quad x >0,\ \alpha >0,\eqno(1.2)$$

$$\left(I_{-}^\alpha f\right)(x) = {1\over \Gamma(\alpha)} \int_x^\infty (t-x)^{\alpha-1} f(t) dt,\quad x >0,\ \alpha >0,\eqno(1.3)$$
where $\Gamma(z)$ is the Euler gamma-function \cite{luch}.  The corresponding Riemann-Liouville fractional derivatives associated with (1.2), (1.3) are given by

$$\left(D_{0+}^\alpha f\right) (x)= {1\over \Gamma(n-\alpha)} \left({d\over dx}\right)^n \int_0^x {f(t) \over (x-t)^{\alpha-n+1} }dt,\quad x >0,\  n= [\alpha] +1,\eqno(1.4)$$

$$\left(D_{-}^\alpha f\right) (x)= {(-1)^n\over \Gamma(n-\alpha)}  \left({d\over dx}\right)^n \int_x^\infty {f(t) \over (t-x)^{\alpha-n+1}} dt,\quad x >0,\   n= [\alpha] +1,\eqno(1.5)$$
where $[\alpha]$ denotes the integer part of $\alpha >0$.  Our approach will involve the Mellin transform in $L_{\nu, p}(\mathbb{R}_+),\ \nu \in \mathbb{R},\ 1 \le  p \le 2$ (see details in \cite{tit}), which is defined by the integral  
$$f^*(s)= \int_0^\infty f(x) x^{s-1} dx.\eqno(1.6)$$
 It converges in mean with respect to the norm in $L_{p^\prime} (\nu- i\infty, \nu + i\infty),\  p^\prime =p/(p-1)$.   Moreover, the  Parseval equality holds for $f \in L_{\nu, p}(\mathbb{R}_+),\  g \in L_{1-\nu, p^\prime}(\mathbb{R}_+)$
$$\int_0^\infty f(x) g(x) dx= {1\over 2\pi i} \int_{\nu- i\infty}^{\nu+i\infty} f^*(s) g^*(1-s) ds.\eqno(1.7)$$
The inverse Mellin transform is given accordingly
 $$f(x)= {1\over 2\pi i}  \int_{\nu- i\infty}^{\nu+i\infty} f^*(s)  x^{-s} ds,\eqno(1.8)$$
where the integral converges in mean with respect to the norm  in   $L_{\nu, p}(\Omega),\ \Omega = \mathbb{R}_+$ and 
$$||f||_{\nu,p} = \left( \int_{\Omega}  |f(x)|^p x^{\nu p-1} dx\right)^{1/p},\quad 1 \le p < \infty,\eqno(1.9)$$

$$||f||_{\nu,\infty} = \hbox{ess sup}_{x\in \Omega} \ x^\nu  |f(x)|.\eqno(1.10)$$
In particular, letting $\nu= 1/p$ we get the usual space $L_p(\mathbb{R}_+; \  dx)$.   As it is proved in \cite{But} the Mellin convolution transform

$$(Kf)(x)\equiv (k*f)(x) =  \int_0^\infty k\left({x\over u}\right) f(u) {du\over u}\eqno(1.11)$$
belongs to $L_{\nu, p}(\mathbb{R}_+),\ 1 \le p \le 2,\ \nu \in \mathbb{R}$ if $f \in L_{\nu, p}(\mathbb{R}_+)$ and $k \in L_{\nu, 1}(\mathbb{R}_+)$. Moreover, the factorization property in terms of the Mellin transform (1.6) takes place

$$(Kf)^*(s)= k^*(s) f^*(s),\quad s= \nu +i\tau,\ \tau \in \mathbb{R}.\eqno(1.12)$$
When $k \in L_{\nu, r}(\mathbb{R}_+),\  1\le r  < \infty$,  the boundedness of the convolution operator (1.11) $K:  L_{\nu, p}(\mathbb{R}_+) \to L_{\nu, q}(\mathbb{R}_+)$, where

$$ 1 \le p \le \infty,\quad   {1\over q} = {1\over p}+ {1\over r}- 1 \ge 0,$$
is known as the Young inequality

$$|| k*f ||_{\nu, q}  \le  || k||_{\nu, r} || f ||_{\nu, p}.\eqno(1.13) $$

Considering the fractional integration and differentiation as reciprocal operations, we will appeal to the space of functions $f$ $AC^n\left(\mathbb{R}_+\right),\ n \in \mathbb{N}$ 
(see \cite{Samka}, Section 1.1) which have continuous derivatives up to order $n-1$ with $f^{(n-1)} \in AC \left(\mathbb{R}_+\right)$, where  
$AC \left(\mathbb{R}_+\right)$ is the space of absolutely continuous functions.  Moreover, denoting by $I_{0+}^\alpha\left(L_1\right),\ \alpha >0$ the space of functions $f$, represented
by the left-sided fractional integral (1.2) of a summable function $f=  I_{0+}^\alpha \varphi,\ \varphi \in L_1\left(\mathbb{R}_+\right)$, we use its description via Theorem 2.3 in  \cite{Samka} to 
say that $f \in I_{0+}^\alpha\left(L_1\right)$ if and only if $I_{0+}^{n-\alpha} f \in  AC^n\left(\mathbb{R}_+\right),\  n= [\alpha] +1$  and $\left(I_{0+}^{n-k-\alpha} f\right)(0)=0,\ k=0,1,\dots, n-1.$

Let $f \in L_{\nu,p}\left(\mathbb{R}_+\right),\ \nu < 1,\  1\le p\le \infty$,  and for positive integers $n$ let us write  the composition of fractional integrals (1.2) $\left(I_{0+}^n x^{-n} I_{0+}^n\right) f$ as follows

$$\left( I_{0+}^n x^{-n} I_{0+}^n\right) f = {1\over [(n-1)!]^2} \int_0^x {(x-t)^{n-1} \over t^n }\int_0^t (t-u)^{n-1} f(u) du dt$$

$$ = {1\over [(n-1)!]^2} \int_0^x  f(u) \int_u^x (x-t)^{n-1} (t-u)^{n-1} t^{-n} dt du,\eqno(1.14)$$
where the interchange of the order of integration for each $x >0$ is allowed via Fubini's theorem due to the estimate

$$\int_0^x {(x-t)^{n-1} \over t^n }\int_0^t (t-u)^{n-1} |f(u)| du dt $$

$$\le  ||f||_{\nu,p} \int_0^x {(x-t)^{n-1} \over t^n } \left( \int_0^t (t-u)^{q(n-1)} u^{(1-\nu)q-1} du\right)^{1/q} dt $$ 

$$=  ||f||_{\nu,p} \int_0^x (x-t)^{n-1} t^{-\nu} dt   \left( \int_0^1 (1-u)^{q(n-1)} u^{(1-\nu)q-1} du\right)^{1/q}  $$ 

$$=   x^{n-\nu}  ||f||_{\nu,p}\   B \left( 1-\nu, n\right) B^{1/q} \left( (1-\nu)q, 1+ q(n- 1)\right)  < \infty,$$ 
where $B(a,b)$ is Euler's beta-function \cite{luch}.  Returning to (1.13), we calculate the inner integral on the right-hand side of the latter equality via Entry 2.2.6.1 in \cite{prud}, Vol. I to obtain

$$ \int_u^x (x-t)^{n-1} (t-u)^{n-1} t^{-n} dt = {\Gamma^2(n)\  (x-u)^{2n-1}\over   \Gamma(2n) \ u^{n}}\  {}_2F_1 \left( n,\ n;\ 2n; \ 1 -{x\over u} \right),\  x > u >0,\eqno(1.15)$$
where ${}_2F_1(a,b,c; z)$ is the Gauss hypergeometric function \cite{prud}, Vol. III.  Hence we derive finally 

$$\left( I_{0+}^n x^{-n} I_{0+}^n\right) f =  {1\over \Gamma(2n)} \int_0^x (x-u)^{2n-1} u^{-n}  {}_2F_1 \left( n,\ n;\ 2n; \ 1 -{x\over u} \right) f(u) du,\ x >0.\eqno(1.16)$$
Analogously, we treat a similar composition, involving fractional integral (1.3). Precisely, it gives

$$\left( I_{-}^n x^{-n} I_{-}^n\right) f = {1\over [(n-1)!]^2} \int_x^\infty {(t-x)^{n-1} \over t^n }\int_t^\infty (u-t)^{n-1} f(u) du dt$$

$$ = {1\over [(n-1)!]^2} \int_x^\infty  f(u) \int_x^u (t-x)^{n-1} (u-t)^{n-1} t^{-n} dt du,\eqno(1.17)$$
where the interchange of the order of integration is guaranteed by the estimate

$$ \int_x^\infty {(t-x)^{n-1} \over t^n }\int_t^\infty (u-t)^{n-1} |f(u)| du dt $$

$$\le ||f||_{\nu,p} \int_x^\infty {(t-x)^{n-1} \over t^n } \left( \int_t^\infty  (u-t)^{q(n-1)} u^{(1-\nu)q-1} du\right)^{1/q} dt $$ 

$$=  ||f||_{\nu,p} \int_x^\infty (t-x)^{n-1} t^{-\nu} dt   \left( \int_1^\infty  (u-1)^{q(n-1)} u^{(1-\nu)q-1} du\right)^{1/q}  $$ 

$$=   x^{n-\nu}  ||f||_{\nu,p}\   B \left( \nu-n, n\right) B^{1/q} \left( (\nu-n)q, 1+ q(n- 1)\right)  < \infty$$ 
when $f \in L_{\nu,p}\left(\mathbb{R}_+\right),\ \nu > n,\  1\le p\le \infty$. Therefore we obtain from (1.17)

$$\left( I_{-}^n x^{-n} I_{-}^n\right) f =  {x^{-n}\over \Gamma(2n)} \int_x^\infty (u-x)^{2n-1}   {}_2F_1 \left( n,\ n;\ 2n; \ 1 -{u\over x} \right) f(u) du,\ x >0.\eqno(1.18)$$
We will call compositions (1.16), (1.18) the Laguerre $n$-fold integration operators, namely,

$$(L_{0+}^n f)(x)=    {1\over \Gamma(2n)} \int_0^x (x-u)^{2n-1} u^{-n}  {}_2F_1 \left( n,\ n;\ 2n; \ 1 -{x\over u} \right) f(u) du,\ x >0,\eqno(1.19)$$

$$(L_{-}^n f)(x)=    {x^{-n}\over \Gamma(2n)} \int_x^\infty (u-x)^{2n-1}   {}_2F_1 \left( n,\ n;\ 2n; \ 1 -{u\over x} \right) f(u) du,\ x >0.\eqno(1.20)$$
Hence,  taking into account (1.1),  there hold the relations

$$\theta^n  L_{0+}^n f = f, \quad\quad   \theta^n  L_{-}^n f = f.\eqno(1.21)$$

The main goal of this paper is to investigate fractional analogs of the Laguerre operators (1.19), (1.20) when the positive integer $n$ is replaced by the real positive $\alpha$. We will study their mapping properties in spaces $L_{\nu,p}\left(\mathbb{R}_+\right)$, establish semigroup properties and formulas of the integration by parts as well as their Mellin-Barnes representations. Finally we will apply these results to the solvability of the corresponding second kind integral equations.

\section{Laguerre fractional integrals and  their properties} 

Basing on definitions (1.19), (1.20) of the Laguerre integration operators for positive integers, we define their fractional counterparts on the positive half-axis, respectively, as follows

$$(L_{0+}^\alpha f)(x)=    {1\over \Gamma(2\alpha)} \int_0^x (x-u)^{2\alpha-1} u^{-\alpha}  {}_2F_1 \left( \alpha,\ \alpha;\ 2\alpha; \ 1 -{x\over u} \right) f(u) du,\eqno(2.1)$$

$$(L_{-}^\alpha f)(x)=    {x^{-\alpha}\over \Gamma(2\alpha)} \int_x^\infty (u-x)^{2\alpha-1}   {}_2F_1 \left( \alpha,\ \alpha;\ 2\alpha; \ 1 -{u\over x} \right) f(u) du,\eqno(2.2)$$
where $\alpha >0.$  Furthermore,  fractional integrals (2.1), (2.2)  can be treated as the Mellin convolution transform (1.11) with the Gauss hypergeometric function as the kernel. For instance, writing (2.1) in the form

$$(L_{0+}^\alpha f)(x)=    {1\over \Gamma(2\alpha)} \int_0^\infty  \left({x\over u} - 1\right)_+^{2\alpha-1}   {}_2F_1 \left( \alpha,\ \alpha;\ 2\alpha; \ 1 -{x\over u} \right) f(u) u^{\alpha-1} du,\eqno(2.3)$$
where

$$(x-1)_+^{2\alpha-1} = \begin{cases} (x-1)^{2\alpha-1},\quad x \ge 1,\\ 0 \quad\quad\quad\quad,\quad x < 1,\end{cases}$$
we find that it represents the  Mellin convolution (1.11) of the function $x^\alpha f(x)$ with the kernel $k_+(x)= \left[ \Gamma(2\alpha) \right]^{-1} (x-1)_+^{2\alpha-1} {}_2F_1 \left( \alpha,\ \alpha;\ 2\alpha; \ 1 - x\right).$  The latter function has the Mellin-Barnes integral representation (see Entry 8.4.49.24 in \cite{prud}, Vol. III), and we  obtain

$${ (x-1)_+^{2\alpha-1} \over \Gamma(2\alpha)} {}_2F_1 \left( \alpha,\ \alpha;\ 2\alpha; \ 1 - x\right) = {1\over 2\pi i} \int_{\gamma-i\infty}^{\gamma+i\infty} \left( {\Gamma(1-\alpha-s)\over \Gamma(1-s)}\right)^2 x^{-s} ds, \eqno(2.4)$$
where $x >0,\  \alpha >0,\  \gamma < 1- \alpha.$ Analogously,  fractional integral (2.2) is the Mellin convolution of the function  $x^\alpha f(x)$ with the kernel $k_-(x)= \left[ \Gamma(2\alpha) \right]^{-1}  x^{-\alpha} (1-x)_+^{2\alpha-1} {}_2F_1 \left( \alpha,\ \alpha;\ 2\alpha; \ 1 - x^{-1} \right).$  The latter function has the Mellin-Barnes integral representation (cf.  Entry 8.4.49.25 in \cite{prud}, Vol. III), and we  find 

$${ x^{-\alpha} (1-x)_+^{2\alpha-1} \over \Gamma(2\alpha)} {}_2F_1 \left( \alpha,\ \alpha;\ 2\alpha; \ 1 - {1\over x}\right) = {1\over 2\pi i} \int_{\gamma-i\infty}^{\gamma+i\infty} \left( {\Gamma(s)\over \Gamma(s+\alpha)}\right)^2 x^{-s} ds, \eqno(2.5)$$
where $x >0,\  \alpha >0,\  \gamma > 0.$  On the other hand, appealing to the Boltz formula for the Gauss hypergeometric function (see \cite{prud}, Vol. III, Entry 7.3.1.3) and Entry 7.3.1.70, these kernels can be written in terms of the associated Legendre functions $P_\nu^\mu(z)$ \cite{prud}, Vol. III due to the identity 

$$   {x^{-\alpha}  (x-1)^{2\alpha-1} \over \Gamma(2\alpha)} {}_2F_1 \left( \alpha,\ \alpha;\ 2\alpha; \ 1 - {1\over x} \right) $$

$$=  {\sqrt\pi   \over  x^{1/4} \  \Gamma(\alpha)} \   (x-1)^{\alpha-1/2}   P_{-1/2}^{1/2-\alpha} \left({1+x\over 2\sqrt x}\right), \quad x > 1.\eqno(2.6)$$
Taking into account asymptotic behavior of the Gauss hypergeometric function at infinity \cite{prud}, Vol. III, we see that kernels of the Laguerre fractional integrals (2.1), (2.2) behave as follows

$$  k_+(x)=  { (x-1)_+^{2\alpha-1} \over \Gamma(2\alpha)} {}_2F_1 \left( \alpha,\ \alpha;\ 2\alpha; \ 1 - x\right) = O\left( x^{\alpha-1} \log(x)\right),\ x \to +\infty,\eqno(2.7)$$

$$k_-(x)= { x^{-\alpha} (1-x)_+^{2\alpha-1} \over \Gamma(2\alpha)} {}_2F_1 \left( \alpha,\ \alpha;\ 2\alpha; \ 1 - {1\over x}\right) = O\left( \log(x)\right),\ x \to 0+.\eqno(2.8)$$
Therefore $k_+(x) \in L_{\nu,1}\left(\mathbb{R}_+\right)$ when $\alpha >0,\ \alpha +\nu < 1$ and  $k_-(x) \in L_{\nu,1}\left(\mathbb{R}_+\right)$ when, in turn,  $\alpha >0,\  \nu > 0$.  Hence if $f \in L_{\alpha+\nu,p}\left(\mathbb{R}_+\right),\ 1\le p\le \infty$ we apply the generalized Minkowski inequality and take into account from (1.15) that the hypergeometric function  ${}_2F_1 \left( \alpha,\ \alpha;\ 2\alpha; \ 1 - x\right),\ \alpha > 0, \ x \ge 1$ is positive to obtain

$$\left|\left| L_{0+}^\alpha f \right|\right|_{\nu,p} = \left( \int_0^\infty \left| {1\over \Gamma(2\alpha)} \int_0^\infty  \left(u - 1\right)_+^{2\alpha-1}   {}_2F_1 \left( \alpha,\ \alpha;\ 2\alpha; \ 1 - u \right) f\left({x\over u}\right) \left({x\over u}\right)^{\alpha} {du\over u} \right|^p\right.$$

$$\left. \times  x^{\nu p -1} dx \right)^{1/p} \le   {1\over \Gamma(2\alpha)} \int_0^\infty  \left(u - 1\right)_+^{2\alpha-1}  u^{-\alpha}  {}_2F_1 \left( \alpha,\ \alpha;\ 2\alpha; \ 1 - u \right)  $$

$$\times \left( \int_0^\infty  \left| f\left({x\over u}\right) \right|^p x^{(\alpha+\nu) p -1} dx \right)^{1/p}  {du\over u} =  C_+  \left|\left|  f \right|\right|_{\alpha+\nu,p},\eqno(2.9)$$
where

$$C_+ =   {1\over \Gamma(2\alpha)} \int_1^\infty  \left(u - 1\right)^{2\alpha-1}  u^{\nu-1}  {}_2F_1 \left( \alpha,\ \alpha;\ 2\alpha; \ 1 - u \right)  du,\ \alpha >0,\ \alpha+\nu < 1.\eqno(2.10)$$
Analogously,

$$\left|\left| L_{-}^\alpha f \right|\right|_{\nu,p} = \left( \int_0^\infty \left| {1\over \Gamma(2\alpha)} \int_0^\infty  \left(1-u\right)_+^{2\alpha-1} u^{-\alpha }  {}_2F_1 \left( \alpha,\ \alpha;\ 2\alpha; \ 1 - {1\over u} \right) f\left({x\over u}\right) \left({x\over u}\right)^{\alpha} {du\over u} \right|^p \right.$$

$$\left. \times x^{\nu p -1} dx \right)^{1/p} \le   {1\over \Gamma(2\alpha)} \int_0^\infty  \left(1-u\right)_+^{2\alpha-1}  u^{- 2\alpha}  {}_2F_1 \left( \alpha,\ \alpha;\ 2\alpha; \ 1 - {1\over u} \right) $$

$$\times \left( \int_0^\infty  \left| f\left({x\over u}\right) \right|^p x^{(\alpha+\nu) p -1} dx \right)^{1/p}  {du\over u} =  C_-  \left|\left|  f \right|\right|_{\alpha+\nu,p},\eqno(2.11)$$
where

$$C_- =   {1\over \Gamma(2\alpha)} \int_0^1  \left(1-u \right)^{2\alpha-1}  u^{\nu-\alpha -1}  {}_2F_1 \left( \alpha,\ \alpha;\ 2\alpha; \ 1 - {1\over u} \right)  du,\ \alpha >0,\  \nu > 0.\eqno(2.12)$$
Furthermore, it is easily seen via (2.7), (2.8) that the kernels $k_+, k_-$ belong to $L_{\nu, r}(\mathbb{R}_+),\  1\le r  < \infty$ under the same conditions. Therefore, appealing to (1.13), we get the following inequalities  for the Laguerre fractional integrals (2.1), (2.2) 

$$ \left|\left| L_{0+}^\alpha f \right|\right|_{\nu,q} \le   || k_+||_{\nu, r} || f ||_{\alpha+\nu, p},\eqno(2.13)$$

$$ \left|\left| L_{-}^\alpha f \right|\right|_{\nu,q} \le   || k_-||_{\nu, r} || f ||_{\alpha+\nu, p},\eqno(2.14)$$
where $ 1 \le p \le \infty,\quad    q^{-1} =  p^{-1} +  r^{-1} - 1 \ge 0.$   These results can be used to establish an analog of the formula for fractional integration by parts, involving  operators (2.1), (2.2)

$$\int_0^\infty f(x) (L_{0+}^\alpha g)(x) dx = \int_0^\infty g(x) (L_{-}^\alpha f)(x) dx.\eqno(2.15)$$
In fact, for  smooth  functions $f, g$ with compact support on $\mathbb{R}_+$  relation (2.15) is verified directly by substitution (2.1) or (2.2) and changing  the order of integration.
Generally, let $f \in   L_{\alpha+\nu, p}(\mathbb{R}_+),\ g \in   L_{1-\nu, r}(\mathbb{R}_+),\  \alpha >0, \ \nu > 0,\  1\le p\le \infty,\ 1 \le r \le \infty,\  p^{-1}+ r^{-1} \ge 1$.  Hence, the H{\"o}lder inequality yields 

$$\left|\int_0^\infty f(x) (L_{0+}^\alpha g)(x) dx \right| \le  || f ||_{\alpha+\nu, p}  \left|\left| L_{0+}^\alpha g \right|\right|_{1-\alpha-\nu,p^\prime}.\eqno(2.16)$$
But since $p^{-1}+ r^{-1} \ge 1$, then $p^\prime \ge r$, and hence (see (2.13))

$$\left|\left| L_{0+}^\alpha g \right|\right|_{1-\alpha-\nu,p^\prime} \le K_1 || g ||_{1-\nu, r},\eqno(2.17)$$
where $K_1=  || k_+||_{1-\alpha-\nu, q}$ and $ q \ge 1$ is defined by the equality $q^{-1}= [ p^\prime]^{-1} - r^{-1} +1.$ Substituting this estimate into (2.16), we get

$$\left|\int_0^\infty f(x) (L_{0+}^\alpha g)(x) dx \right| \le  K_1  || f ||_{\alpha+\nu, p}  || g ||_{1-\nu, r}.\eqno(2.18)$$
It means that the left-hand side of (2.15) represents a bounded bilinear functional on $  L_{\alpha+\nu, p}(\mathbb{R}_+) \times  L_{1-\nu, r}(\mathbb{R}_+).$ On the other hand, the right-hand side of (2.15) can be treated accordingly

$$\left|  \int_0^\infty g(x) (L_{-}^\alpha f)(x) dx\right| \le  || g ||_{1-\nu, r}  \left|\left| L_{-}^\alpha f \right|\right|_{\nu,r^\prime},\quad r^\prime = {r\over r-1}.$$
Then similar to (2.17) we have 

$$ \left|\left| L_{-}^\alpha f \right|\right|_{\nu,r^\prime} \le K_2   || f ||_{\alpha+\nu, p},$$
where  $K_2=  || k_-||_{\nu, q}$ and $ q \ge 1$ is defined by the equality $q^{-1}= [ r^\prime]^{-1} - p^{-1} +1.$ Hence

$$ \left|  \int_0^\infty g(x) (L_{-}^\alpha f)(x) dx\right| \le  K_2   || f ||_{\alpha+\nu, p}  || g ||_{1-\nu, r},$$
and the right-hand side of (2.17) is a bounded bilinear functional on $  L_{\alpha+\nu, p}(\mathbb{R}_+) \times  L_{1-\nu, r}(\mathbb{R}_+).$  Consequently, we proved the following

{\bf Theorem 1.} {\it  Let $\alpha, \nu  > 0, \ 1\le p\le \infty,\   1\le r\le \infty $ be such that $ p^{-1}+ r^{-1} \ge 1$. If $ f \in L_{\alpha+\nu, p}(\mathbb{R}_+),\ g \in 
L_{1-\nu, r}(\mathbb{R}_+)$, then there holds formula $(2.15)$ of fractional integration by parts for Laguerre's fractional integrals $(2.1), (2.2)$.} 

Further, recalling (1.12), (2.4), (2.5), we write the corresponding relations for the Mellin transform of the fractional integrals (2.1), (2.2), namely,

$$(L_{0+}^\alpha f)^*(s) = \left( {\Gamma(1-\alpha-s)\over \Gamma(1-s)}\right)^2  f^*(s+\alpha),\ s= \nu +i\tau,\ \tau \in \mathbb{R},\eqno(2.19)$$

$$(L_{-}^\alpha f)^*(s) = \left( {\Gamma(s)\over \Gamma(s+\alpha)}\right)^2  f^*(s+\alpha),\ s= \nu +i\tau,\ \tau \in \mathbb{R}\eqno(2.20)$$
under the conditions $f \in L_{\alpha+\nu,p}\left(\mathbb{R}_+\right),\ 1\le p\le 2,\ \alpha, \nu >0,\ \alpha+\nu < 1.$  However, when $ \alpha > 1/(2p)$, the Stirling asymptotic formula suggests that the square of the quotient of the  gamma-functions belongs to $L_p(\nu-i\infty, \nu+i\infty),\  p \ge 1.$ Then from (2.19), (2.20) and (1.8) we derive, reciprocally, 

$$(L_{0+}^\alpha f)(x) = {1\over 2\pi i } \int_{\nu-i\infty}^{\nu + i\infty} \left( {\Gamma(1-\alpha-s)\over \Gamma(1-s)}\right)^2  f^*(s+\alpha) x^{-s} ds,\eqno(2.21)$$

$$(L_{-}^\alpha f)(x) =  {1\over 2\pi i } \int_{\nu-i\infty}^{\nu + i\infty} \left( {\Gamma(s)\over \Gamma(s+\alpha)}\right)^2  f^*(s+\alpha) x^{-s} ds,\eqno(2.22)$$
where integrals (2.21), (2.22) converge absolutely due to the H{\"o}lder inequality and mapping $L_p$-properties of the Mellin transform (1.6) (see above).  These formulas are key ingredients to establish the semigroup property for the Laguerre fractional integrals (2.1), (2.2)

$$ L_{0+}^{\alpha+\beta} f = L_{0+}^\alpha L_{0+}^\beta f,\eqno(2.23)$$

$$ L_{-}^{\alpha+\beta} f = L_{-}^\alpha L_{-}^\beta f.\eqno(2.24)$$
Indeed,  the right-hand side of (2.24) can be treated,  substituting  the Mellin-type representation (2.21) for  $(L_{0+}^\beta  f )(x)$ into (2.1) and changing the order of integration by Fubini's theorem owing to the following estimate for each $x >0$

$$ \int_0^x (x-u)^{2\alpha-1} u^{-\alpha-\nu}  {}_2F_1 \left( \alpha,\ \alpha;\ 2\alpha; \ 1 -{x\over u} \right)  \int_{\nu-i\infty}^{\nu + i\infty} \left| \left({\Gamma(1-\beta-s)\over \Gamma(1-s)}\right)^2  f^*(s+\beta) ds \right|  du$$

$$= x^{\alpha-\nu}  \int_0^\infty (y-1)_+^{2\alpha-1} y^{\nu-\alpha-1}  {}_2F_1 \left( \alpha,\ \alpha;\ 2\alpha; \ 1 - y \right)  dy $$

$$\times \int_{\nu-i\infty}^{\nu + i\infty} \left| \left( {\Gamma(1-\beta-s)\over \Gamma(1-s)}\right)^2  f^*(s+\beta) ds \right| < \infty,$$ 
when $\nu < 1-\beta,\ \alpha >0, \ \beta >  1/(2p),\ 1 \le p \le 2,\  f \in  L_{\beta+\nu,p}\left(\mathbb{R}_+\right).$  Hence, using (2.4),  we obtain

$$\left(L_{0+}^\alpha L_{0+}^\beta f\right) (x)= {1\over 2\pi i} \int_{\nu-i\infty}^{\nu + i\infty}  \left( {\Gamma(1-\beta-s)\over \Gamma(1+\alpha -s)}\right)^2  f^*(s+\beta)\  x^{\alpha-s} ds,\ x >0.\eqno(2.25)$$ 
But the right-hand side of (2.25) can be expressed in terms of the Mellin convolution (1.11), employing again Entry 8.4.49.24 in \cite{prud}, Vol. III and shift  property of the Mellin transform (1.6). Hence it gives

 $$\left(L_{0+}^\alpha L_{0+}^\beta f\right) (x)= {x^\alpha \over \Gamma(2(\alpha+\beta))} \int_0^\infty  \left({x\over u} - 1\right)_+^{2(\alpha+\beta)-1}   {}_2F_1 \left( \alpha+\beta,\ \alpha+\beta;\ 2(\alpha+\beta); \ 1 - {x\over u} \right) $$
 
 $$\times \left({x\over u}\right)^{-\alpha} f\left(u\right) u^{\beta-1} du = \left( L_{0+}^{\alpha+\beta} f \right)(x),$$
which proves (2.23).  Analogously,  equality (2.24) can be proved, recalling representation (2.22) and Entry 8.4.49.25 in \cite{prud}, Vol. III.  Thus  we get 

{\bf Theorem 2}.  {\it Let $\nu < 1-\beta \ (\nu  > 0),\  \alpha > 0, \   \beta >  1/(2p),\ 1 \le p \le 2, \ f \in  L_{\beta+\nu,p}\left(\mathbb{R}_+\right)$.  Then semigroup properties $(2.23)\  ( (2.24) )$ hold for all $x >0$.}

 \section{Laguerre fractional derivatives and their properties}
 
 In this section we will define the left- and right-sided fractional order derivatives associated with the Laguerre fractional integrals (2.1), (2.2).   We will do it similarly to the Riemann-Liouville fractional derivatives (1.4), (1.5), involving the operator $\theta= DxD$  of the ordinary Laguerre derivative (1.1).  Indeed, recalling (2.1), (2.2), (1.21) we have, correspondingly,
 
 $$\left({\mathcal D}_{0+}^\alpha f\right) (x) \equiv \left(D_{0+}^\alpha x^\alpha D_{0+}^\alpha \right) f = \theta^m (L_{0+}^{m-\alpha} f)(x),\  \alpha >0, \ m= [\alpha] +1,\ x >0,\eqno(3.1)$$
 
 $$\left({\mathcal D}_{-}^\alpha f\right) (x) \equiv \left(D_{-}^\alpha x^\alpha D_{-}^\alpha \right) f  = \theta^m (L_{-}^{m-\alpha} f)(x),\  \alpha >0, \ m= [\alpha] +1,\  x >0.\eqno(3.2)$$
 An alternative definition of fractional derivatives (3.1), (3.2) can be given, basing on properties of various types of numbers in combinatorial analysis from the classical discrete to the fractional case.  In particular, the following identities in \cite{Riordan} give operator relations, involving the falling factorial operator $[ xD]_n = xD (xD-1) \dots (xD -n+1)$ and the operator $x^nD^n$, namely,

 $$ [ xD]_n = \sum_{k=0}^n s(n,k) (xD)^k = x^n D^n,\eqno(3.3)$$
 where $s(n,k)$ are the Stirling numbers of the first kind. Therefore (1.21) suggests the identity
 
 $$\theta^n f = (D^n x^n D^n) f = x^{-n} \left(  \sum_{k=0}^n s(n,k) (xD)^k \right)^2 f.\eqno(3.4)$$
Hence in order to define the fractional power of the operator $\theta$, we will appeal to the so-called Stirling function $s(\alpha,k)$ of the first kind (cf.  \cite{ButKil}) 
 
$$\left. s(\alpha,k) = {1\over k!} {d^k [ u]_\alpha \over du^k}  \right |_{u=0},\  \alpha > 0, k \in \mathbb{N}_0\eqno(3.5)$$
and the falling factorial function $ [ u]_\alpha $ is naturally defined by the formula

$$  [ u]_\alpha = {\Gamma(u+1)\over \Gamma(u+1-\alpha)}.\eqno(3.6)$$
Since $s(n,k) = 0,\ k \ge n+1$ \cite{ButKil}, the fractional  derivative (3.1) can be given from (3.4) in the form

$$\left({\mathcal D}_{0+}^\alpha f\right) (x) = x^{-\alpha} \left(  \sum_{k=0}^\infty  s(\alpha,k) (xD)^k \right)^2 f.\eqno(3.7)$$
Appealing to the Cauchy product for the series, we rewrite (3.7) as follows

$$ \left({\mathcal D}_{0+}^\alpha f\right) (x) = x^{-\alpha}  \sum_{k=0}^\infty  c_k(\alpha) (xD)^k f, \eqno(3.8)$$
where

$$c_k(\alpha) = \sum_{j=0}^k  s(\alpha,j)\  s(\alpha,k-j).\eqno(3.9)$$
In fact, employing the Leibniz differentiation formula for the product of functions, we have from (3.5)

$$\left. c_k(\alpha) =  \lim_{u\to 0} {1\over k!} \sum_{j=0}^k \binom{k}{j} {d^j [ u]_\alpha \over du^j}\  {d^{k-j} [ u]_\alpha \over du^{k-j}} = {1\over k!} {d^k [ u]^2_\alpha \over du^k}  \right |_{u=0} .$$
Therefore we obtain formally from (3.7) (cf. (3.3))

$$ \left. x^{-\alpha}  \sum_{k=0}^\infty  c_k(\alpha) (xD)^k f = x^{-\alpha}  \sum_{k=0}^\infty  {1\over k!} {d^k [ u]^2_\alpha \over du^k}  \right |_{u=0}  (xD)^k f$$

$$=  x^{-\alpha}  \left(  [ xD]_\alpha \right)^2 f  = \left(D_{0+}^\alpha x^\alpha D_{0+}^\alpha \right) f.\eqno(3.10)$$
In the same manner fractional derivative (3.2) can be interpreted.  A more straightforward interpretation can be realized, writing the identity

$$(D^n x^n D^n) f = n! \left(\sum_{k=0}^n \binom{n}{k} {x^k\over k!} \ D^{k+n}\right) f .\eqno(3.11)$$
Hence we set, for instance, 

$$ \left({\mathcal D}_{0+}^\alpha f\right) (x) = \Gamma(\alpha+1)  \left(\sum_{k=0}^\infty { [\alpha]_k \over (k! )^2} x^k \ D_{0+}^{k+\alpha}\right) f .\eqno(3.12)$$
The right-hand side of the latter equality (3.12) can be formally interpreted  in terms of the Gauss hypergeometric function at the unity, and we deduce the following operational relation for the Laguerre fractional derivative (3.1) in terms of the Riemann-Louville derivative (1.4)

$$ \left({\mathcal D}_{0+}^\alpha f\right) (x) =  {\Gamma(1+\alpha+ xD)\over \Gamma(1+ xD)}\  D_{0+}^{\alpha} f .\eqno(3.13)$$
Analogously,  for the derivative (3.2) we find 

$$ \left({\mathcal D}_{-}^\alpha f\right) (x) =  {\Gamma(1+\alpha+ xD)\over \Gamma(1+ xD)}\  D_{-}^{\alpha} f .\eqno(3.14)$$
Now,  recalling (2.21), (2.22), we write fractional integrals in (3.1), (3.2)  accordingly, 

$$(L_{0+}^{m-\alpha} f)(x) =  {1\over 2\pi i } \int_{\nu-i\infty}^{\nu + i\infty} \left( {\Gamma(1-m+\alpha -s)\over \Gamma(1-s)}\right)^2  f^*(s+ m-\alpha) x^{-s} ds,\eqno(3.15)$$ 

$$(L_{-}^{m-\alpha} f)(x) = {1\over 2\pi i } \int_{\gamma-i\infty}^{\gamma + i\infty} \left( {\Gamma(s)\over \Gamma(s+ m-\alpha)}\right)^2  f^*(s+m-\alpha) x^{-s} ds.\eqno(3.16)$$
Basing on relations (1.21) and Theorem 2,  it is not difficult to establish the following identities

$$ {\mathcal D}_{0+}^\alpha L_{0+}^\alpha f  = f\eqno(3.17)$$
under conditions  $  \nu < 1-\alpha,\ \alpha > 1/(2p),  \ f \in  L_{\alpha+\nu,p}\left(\mathbb{R}_+\right) \cap L_{\nu,p}\left(\mathbb{R}_+\right),\ m= [\alpha] +1,\ 1 \le p \le 2,$

$$ {\mathcal D}_{-}^\alpha L_{-}^\alpha f  = f\eqno(3.18)$$
when  $  \alpha > 0,  \ f \in  L_{\alpha+\gamma,p}\left(\mathbb{R}_+\right) \cap L_{\gamma,p}\left(\mathbb{R}_+\right),\ \gamma > m,\ m= [\alpha] +1,\ 1 \le p \le 2.$ Let  $s^{2\alpha} f^*(s+m-\alpha) \in L_1\left(\nu-i\infty,\ \nu+i\infty\right) \cap  L_1\left(\gamma-i\infty,\ \gamma+i\infty\right)$.  Hence via the asymptotic of the quotient of gamma-functions and the Lebesgue dominated convergence theorem we find that   integrals (3.15), (3.16)  converge absolutely.   Then, applying to the both sides of (3.15), (3.16)   the differential operator $\theta^m$, we differentiate under the integral sign on their right-hand sides due to the absolute and uniform convergence by $x \ge x_0 >0$ to obtain, involving (3.1), (3.2),  the equalities

$$  \left({\mathcal D}_{0+}^\alpha f\right) (x) =  {1\over 2\pi i } \int_{\nu-i\infty}^{\nu + i\infty} \left( {\Gamma(1+m-\alpha-s) (s)_m \over \Gamma(1-s)}  \right)^2  f^*(s+m-\alpha) x^{-s-m} ds,$$ 

$$ \left({\mathcal D}_{-}^\alpha f\right) (x) = {1\over 2\pi i } \int_{\gamma-i\infty}^{\gamma + i\infty} \left( {\Gamma(s) (s)_m \over \Gamma(s+ m-\alpha)}\right)^2  f^*(s+m-\alpha) x^{-s-m} ds,$$
where $(s)_m = s(s+1)\dots (s+m-1)$ is the Pochhammer symbol.  It can be simplified, employing the addition and reflection formulas for the gamma-function with a simple substitution.   In fact, we derive

$$  \left({\mathcal D}_{0+}^\alpha f\right) (x) =  {1\over 2\pi i } \int_{\nu-i\infty}^{\nu + i\infty} \left( {\Gamma(1+m-\alpha -s) \over \Gamma(1-s-m)}  \right)^2  f^*(s+m-\alpha) x^{-s-m} ds$$

$$=  {1\over 2\pi i } \int_{\nu+m-i\infty}^{\nu+m + i\infty} \left( {\Gamma(1+\alpha-s) \over \Gamma(1-s)}  \right)^2  f^*(s-\alpha) x^{-s} ds.\eqno(3.19)$$ 
Analogously,  we find

$$ \left({\mathcal D}_{-}^\alpha f\right) (x) = {1\over 2\pi i } \int_{\gamma+m-i\infty}^{\gamma +m+ i\infty} \left( {\Gamma(s )  \over \Gamma(s -\alpha)}\right)^2  f^*(s-\alpha) x^{-s} ds.\eqno(3.20)$$
The differentiation under the integral sign in (3.19), (3.20) is motivated by the estimates, respectively, 

$$ \int_{\nu-i\infty}^{\nu + i\infty} \left| {\Gamma(1+m-\alpha -s) \over \Gamma(1-s-m)}  \right|^2  \left| f^*(s+m-\alpha) x^{-s-m} ds\right| $$

$$\le C_m\  x_0^{-\nu-m} \int_{\nu-i\infty}^{\nu + i\infty}|s|^{2\alpha}  \left| f^*(s+m-\alpha)  ds\right| < \infty,\eqno(3.21)$$

$$\int_{\gamma-i\infty}^{\gamma + i\infty} \left| {\Gamma(s +m)  \over \Gamma(s +m-\alpha)}\right|^2 \left| f^*(s+m-\alpha) x^{-s-m} ds\right|$$

$$\le C_m\  x_0^{-\nu-m} \int_{\gamma-i\infty}^{\gamma + i\infty}|s|^{2\alpha}  \left| f^*(s+m-\alpha)  ds\right| < \infty,\eqno(3.22)$$
where $C_m >0 $ is an absolute constant.  Moreover, comparing with (2.21), (2.22), we identify Laguerre fractional derivatives in terms of the corresponding fractional integrals. Precisely, we have

$$ \left({\mathcal D}_{0+}^\alpha f\right) (x) \equiv  (L_{0+}^{-\alpha} f)(x),\quad\quad  \left({\mathcal D}_{-}^\alpha f\right) (x) \equiv  (L_{-}^{-\alpha} f)(x).\eqno(3.23)$$
Furthermore, reciprocally to (3.17), (3.18), one can prove, recalling (2.1), (2.2),  (3.19), (3.20) and Fubini's theorem, the identities

$$ L_{0+}^\alpha {\mathcal D}_{0+}^\alpha  f  = f,\quad\quad\quad   L_{-}^\alpha {\mathcal D}_{-}^\alpha  f  = f.\eqno(3.24)$$
Indeed, we find from (2.4) and (1.8)

$$ \left(L_{0+}^\alpha {\mathcal D}_{0+}^\alpha  f\right) (x)  = {1\over \Gamma(2\alpha)} \int_0^x (x-u)^{2\alpha-1} u^{-\alpha}  {}_2F_1 \left( \alpha,\ \alpha;\ 2\alpha; \ 1 -{x\over u} \right) $$

$$\times  {1\over 2\pi i } \int_{\nu+m-i\infty}^{\nu+m + i\infty} \left( {\Gamma(1+\alpha-s) \over \Gamma(1-s)}  \right)^2  f^*(s-\alpha) u^{-s} ds du$$

$$=  {1\over 2\pi i } \int_{\nu+m-i\infty}^{\nu+m + i\infty} \left( {\Gamma(1+\alpha-s) \over \Gamma(1-s)}  \right)^2  \left( {\Gamma(1-s) \over \Gamma(1+\alpha-s)}  \right)^2 f^*(s-\alpha) x^{\alpha-s} ds $$

$$=  {1\over 2\pi i } \int_{\nu+m-i\infty}^{\nu+m + i\infty}  f^*(s-\alpha) x^{\alpha-s} ds =  {1\over 2\pi i } \int_{\nu+m-\alpha -i\infty}^{\nu+m -\alpha + i\infty}  f^*(s) x^{-s} ds = f(x),$$ 
where the interchange of the order of integration is allowed by virtue of the estimate (see (2.9))

$$\int_0^1 (1-u)^{2\alpha-1} u^{-\alpha-\nu-m}  {}_2F_1 \left( \alpha,\ \alpha;\ 2\alpha; \ 1 -{1\over u} \right) $$

$$\times   \int_{\nu+m-i\infty}^{\nu+m + i\infty} \left|\left( {\Gamma(1+\alpha-s) \over \Gamma(1-s)}  \right)^2  f^*(s-\alpha)  ds\right|  du < \infty$$
under conditions  $f \in L_{\nu+m-\alpha, 1}\left(\mathbb{R}_+\right),\ s^{2\alpha} f^*(s+m-\alpha) \in L_1\left(\nu-i\infty,\ \nu+i\infty\right),\ m= [\alpha]+1,\  \alpha >0, \ \nu < 1-m.$ In the same manner we establish the second identity in (3.24), assuming that  $f \in L_{\gamma+m-\alpha, 1}\left(\mathbb{R}_+\right),\ s^{2\alpha} f^*(s+m-\alpha) \in L_1\left(\gamma-i\infty,\ \gamma+i\infty\right),\ m= [\alpha]+1,\  \alpha >0, \ \gamma> \alpha-m.$

An analog of formula (2.15) for the Laguerre fractional derivatives can be established by the following theorem. 

{\bf Theorem 3}.  {\it Let $\alpha >0,\  \nu < 1+ \alpha-m, \ m =[\alpha]+1,   f  \in L_{1-\nu-m,1}\left(\mathbb{R}_+\right), g \in   L_{\nu+m- \alpha,1}\left(\mathbb{R}_+\right)$,\    $s^{2\alpha} g^*(s-\alpha) \in L_1\left(\nu+m -i\infty, \nu+m + i\infty\right),\   s^{2\alpha} f^*(s-  \alpha) \in   L_1\left(1+\alpha-\nu -m-i\infty, 1+\alpha-\nu-m+ i\infty\right)$.  Then the following identity holds}

$$\int_0^\infty f(x) ({\mathcal D}_{0+}^\alpha g)(x) dx = \int_0^\infty g(x) ({\mathcal D}_{-}^\alpha f)(x) dx.\eqno(3.25)$$
\begin{proof}  Indeed,  taking  representation (3.19) for the Laguerre fractional derivative $({\mathcal D}_{0+}^\alpha g)(x)$ and plugging it  into  the left-hand side of (3.25), we change the order of integration by Fubini's theorem due to the involved assumptions and estimates (3.21), (3.22).   Hence we obtain via (1.6), (1.7), (3.20)   and a simple substitution 

$$\int_0^\infty f(x) ({\mathcal D}_{0+}^\alpha g)(x) dx =  
{1\over 2\pi i } \int_{\nu+m-i\infty}^{\nu+m + i\infty}   \left( {\Gamma(1+\alpha-s) \over \Gamma(1-s)}  \right)^2  g^*(s-\alpha)  f^*(1-s) ds$$

$$= {1\over 2\pi i } \int_{\gamma+m-i\infty}^{\gamma+m+ i\infty}   \left( {\Gamma(s) \over \Gamma(s-\alpha)}  \right)^2  g^*(1-s)  f^*(s-\alpha) ds 
=  \int_0^\infty g(x) ({\mathcal D}_{-}^\alpha f)(x) dx,$$
where $\gamma =  1+\alpha- 2m-\nu.$

\end{proof} 

\section{The Volterra-type integral equation of the second kind }

Let us consider the following Volterra-type equation of the second kind, involving the Laguerre fractional integral (2.1) in the space $L_{\nu,p} (0,l),\ 1\le p \le \infty,\  l > 0$

$$f(x)= g(x) +  {\lambda \over \Gamma(2\alpha)} \int_0^x (x-u)^{2\alpha-1} u^{-\alpha}  {}_2F_1 \left( \alpha,\ \alpha;\ 2\alpha; \ 1 -{x\over u} \right) f(u) du,\quad x \in (0,l),\eqno(4.1)$$
where the parameter $\lambda \in \mathbb{C}$, \ $g \in L_{\nu,p} (0,l)$ is a given function and $f$ is to  be determined in the space $L_{\nu,p} (0,l)$.   In order to study the solvability of the equation (4.1) in this space we first show that the integral operator $L_{0+}^{\alpha} : L_{\nu,p} (0,l) \to L_{\nu,p} (0,l)$ is bounded.  In fact, we have similarly to (2.9)

$$\left|\left| L_{0+}^\alpha f \right|\right|_{\nu,p} = \left( \int_0^l \left| {1\over \Gamma(2\alpha)} \int_0^x  \left({x\over u} - 1\right)^{2\alpha-1}   {}_2F_1 \left( \alpha,\ \alpha;\ 2\alpha; \ 1 - {x\over u} \right) f\left(u\right)  u^{\alpha-1}  du \right|^p  x^{\nu p -1} dx \right)^{1/p} $$

$$= \left( \int_0^l \left| {1\over \Gamma(2\alpha)} \int_1^\infty  \left(u - 1\right)^{2\alpha-1}   {}_2F_1 \left( \alpha,\ \alpha;\ 2\alpha; \ 1 - u  \right) f\left({x\over u}\right)  \left({x\over u}\right)^{\alpha}  {du \over u}\right|^p  x^{\nu p -1} dx \right)^{1/p} $$

$$\le   {1\over \Gamma(2\alpha)} \int_1^\infty  \left(u - 1\right)^{2\alpha-1}  u^{\nu-1} {}_2F_1 \left( \alpha,\ \alpha;\ 2\alpha; \ 1 - u \right)   \left( \int_0^{l/u}   \left| f\left(x\right) \right|^p x^{(\alpha+\nu) p -1} dx \right)^{1/p}  du$$

$$\le  C_+ \  l^\alpha \left|\left|  f \right|\right|_{\nu,p},\eqno(4.2)$$
where $\alpha >0,\ \alpha+\nu < 1$ and $C_+$ is defined by (2.10).   This estimate proves the boundedness of the Laguerre integral (2.1) in the space  $L_{\nu,p} (0,l)$.  Further, equation (4.2) reads in the operator form

$$\left(I - \lambda L_{0+}^\alpha \right) f(x)= g(x),\quad  x \in (0, l).\eqno(4.3)$$
Hence it can be solved in terms  of the Neumann series as follows

$$ f(x)= \left(I - \lambda L_{0+}^\alpha \right)^{-1} g(x) = \left( 1+ \sum_{n=1}^\infty \lambda^n \left(  L_{0+}^\alpha \right)^n \right) g(x),\eqno(4.4)$$
whose absolute convergence is guaranteed owing to (4.2) in the open disk $| \lambda| < (C_+ l^\alpha)^{-1}$ since

$$\left|\left|  1+ \sum_{n=1}^\infty \lambda^n \left(  L_{0+}^\alpha \right)^n \right|\right|_{\nu,p} \le   \sum_{n=0}^\infty |\lambda|^n  \left| \left|L_{0+}^\alpha \right|\right|^n_{\nu,p}\le 
 \sum_{n=0}^\infty |\lambda|^n  (C_+ l^\alpha)^n $$

$$= {1\over 1 -  |\lambda|\ C_+ l^\alpha },\quad  | \lambda| < (C_+ l^\alpha)^{-1}.$$
This  unique solution can be written in the explicit form if we establish the semigroup property (2.24) in the space $L_{\nu,p} (0,l)$.  To do this, we have

$$\left(L_{0+}^\beta L_{0+}^\alpha f \right) (x)=  {1 \over  \Gamma(2\beta) \Gamma(2\alpha)} \int_0^x (x-u)^{2\beta-1} u^{-\beta}  {}_2F_1 \left( \beta,\ \beta;\ 2\beta; \ 1 -{x\over u} \right) $$

$$\times  \int_0^u (u-t)^{2\alpha-1} t^{-\alpha}  {}_2F_1 \left( \alpha,\ \alpha;\ 2\alpha; \ 1 -{u\over t} \right) f(t) dt du$$

$$=  {1 \over  \Gamma(2\beta) \Gamma(2\alpha)} \int_0^x  t^{-\alpha} f(t) \int_t^x (x-u)^{2\beta-1} (u-t)^{2\alpha-1} u^{-\beta}  {}_2F_1 \left( \beta,\ \beta;\ 2\beta; \ 1 -{x\over u} \right) $$

$$\times  {}_2F_1 \left( \alpha,\ \alpha;\ 2\alpha; \ 1 -{u\over t} \right) du dt,\eqno(4.5)$$
where the interchange of the order of integration is permitted due to the estimates  via the H{\"o}lder and generalized Minkowski inequalities

$$\int_0^x  t^{-\alpha} |f(t) |\int_t^x (x-u)^{2\beta-1} (u-t)^{2\alpha-1} u^{-\beta}  {}_2F_1 \left( \beta,\ \beta;\ 2\beta; \ 1 -{x\over u} \right) $$

$$\times  {}_2F_1 \left( \alpha,\ \alpha;\ 2\alpha; \ 1 -{u\over t} \right) du dt$$

$$ \le ||f||_{\nu,p} \left(\int_0^x  t^{(\alpha-\nu)p^\prime -1} \left(\int_t^x \left({x\over u} - 1\right)^{2\beta-1} \left({u\over t} -1\right)^{2\alpha-1} u^{\beta-1}  {}_2F_1 \left( \beta,\ \beta;\ 2\beta; \ 1 -{x\over u} \right) \right.\right.$$

$$\left.\left.  \times  {}_2F_1 \left( \alpha,\ \alpha;\ 2\alpha; \ 1 -{u\over t} \right) du \right)^{p^\prime} dt\right)^{1/p^\prime}$$

$$ \le ||f||_{\nu,p} \int_0^x  \left(x- u\right)^{2\beta-1} u^{-\beta}  {}_2F_1 \left( \beta,\ \beta;\ 2\beta; \ 1 -{x\over u} \right)  \left( \int_0^u  t^{(1-\alpha-\nu)p^\prime -1}  \left(u-t\right)^{(2\alpha-1)p^\prime}  \right.$$

$$\left. \times  \left( {}_2F_1 \left( \alpha,\ \alpha;\ 2\alpha; \ 1 -{u\over t} \right) \right)^{p^\prime} dt\right)^{1/p^\prime} du$$

$$ =  ||f||_{\nu,p} \int_0^x  \left(x- u\right)^{2\beta-1} u^{\alpha-\nu -\beta}  {}_2F_1 \left( \beta,\ \beta;\ 2\beta; \ 1 -{x\over u} \right)  du $$

$$\times \left( \int_0^1 t^{(1-\alpha-\nu)p^\prime -1}  \left(1-t\right)^{(2\alpha-1)p^\prime}   \left( {}_2F_1 \left( \alpha,\ \alpha;\ 2\alpha; \ 1 -{1\over t} \right) \right)^{p^\prime} dt\right)^{1/p^\prime}  < \infty$$
when $x >0, \  \alpha > 1/(2p),\ \beta > 0, \  \nu < 1,\ 1 \le p \le \infty.$  Hence, returning to (4.5), we calculate the inner integral, employing (2.4), (2.25).  Thus we find

$${1 \over  \Gamma(2\beta) \Gamma(2\alpha)} \int_t^x (x-u)^{2\beta-1} (u-t)^{2\alpha-1} u^{-\beta}  {}_2F_1 \left( \beta,\ \beta;\ 2\beta; \ 1 -{x\over u} \right)$$

$$\times   {}_2F_1 \left( \alpha,\ \alpha;\ 2\alpha; \ 1 -{u\over t} \right) du =   {t^{2\alpha-1}  \over   2\pi i \ \Gamma(2\beta) }\int_0^\infty \left({x\over u}-1\right)_+^{2\beta-1}  u^{\beta-1}  {}_2F_1 \left( \beta,\ \beta;\ 2\beta; \ 1 -{x\over u} \right) $$

$$\times \int_{\gamma-i\infty}^{\gamma+i\infty} \left( {\Gamma(1-\alpha-s)\over \Gamma(1-s)}\right)^2 \left({u\over t}\right)^{-s} ds du = {t^{2\alpha-1}  x^\beta  \over   2\pi i  } \int_{\gamma-i\infty}^{\gamma+i\infty} \left( {\Gamma(1-\alpha-s)\over \Gamma(1+\beta-s)}\right)^2 \left({x\over t}\right)^{-s} ds $$

$$=   { t^{2\alpha+\beta-1} \over \Gamma(2(\alpha+\beta))}  \left({x\over t} - 1\right)_+^{2(\alpha+\beta)-1}   {}_2F_1 \left( \alpha+\beta,\ \alpha+\beta;\ 2(\alpha+\beta); \ 1 - {x\over t} \right),$$
where the interchange of the order of integration is allowed by Fubini's theorem under conditions $ \alpha > 1/2,\ \beta > 0, \ \gamma +\alpha < 1$.  Consequently, substituting this result into the right-hand side of the latter equality (4.5), we obtain finally

$$\left(L_{0+}^\beta L_{0+}^\alpha f \right) (x)=  {1 \over \Gamma(2(\alpha+\beta))}   \int_0^x   \left(x -  t\right)^{2(\alpha+\beta)-1}   t^{-\alpha-\beta}  $$

$$\times  {}_2F_1 \left( \alpha+\beta,\ \alpha+\beta;\ 2(\alpha+\beta); \ 1 - {x\over t} \right) f(t) dt =  \left( L_{0+}^{\alpha+\beta} f \right) (x), $$
completing the proof of the semigroup property (2.24) in the space $L_{\nu,p} (0,l)$.  Therefore, recalling (4.4), we write the unique solution of the integral equation (4.1) in the form

$$f(x)= g(x) + \sum_{n=1}^\infty \lambda^n \left(L_{0+}^{\alpha n} g \right) (x),\quad x \in (0,l),\eqno(4.6)$$
where $g \in L_{\nu,p} (0,l),\ 1 \le p \le \infty,\ \nu < 1-\alpha,\ \alpha > 1/2$ and the series converges absolutely in the open disk $| \lambda| < (C_+ l^\alpha)^{-1}$.  Hence,  using the definition (2.1) of the Laguerre fractional integral, we write the series in (4.6) as follows

$$ \sum_{n=1}^\infty \lambda^n \left(L_{0+}^{\alpha n} g \right) (x) =    \sum_{n=1}^\infty  {\lambda^n \over \Gamma(2\alpha n)} \int_0^x (x-u)^{2\alpha n-1} u^{-\alpha n}  {}_2F_1 \left( \alpha n,\ \alpha n ;\ 2\alpha n; \ 1 -{x\over u} \right) g(u) du $$

$$=   \int_0^x g(u) \sum_{n=1}^\infty  {\lambda^n \over \Gamma(2\alpha n)} \  (x-u)^{2\alpha n-1} u^{-\alpha n}  {}_2F_1 \left( \alpha n,\ \alpha n ;\ 2\alpha n; \ 1 -{x\over u} \right)  du.\eqno(4.7) $$
The term-wise integration and summation is, indeed,  allowed,   and it can be shown,   employing  again the Boltz  formula  for the Gauss hypergeometric function and its representation (2.6) in terms of the associated Legendre function.   This yields the formula 

$${ (x-u)^{2\alpha n-1} \over \Gamma(2\alpha n)} \   u^{-\alpha n}  {}_2F_1 \left( \alpha n,\ \alpha n ;\ 2\alpha n; \ 1 -{x\over u} \right) $$

$$=  {x^{-\alpha n}  (x-u)^{2\alpha n-1} \over \Gamma(2\alpha n)} \   {}_2F_1 \left( \alpha n,\ \alpha n ;\ 2\alpha n; \ 1 -{u\over x} \right) $$

$$ =   {\sqrt\pi  \over (xu)^{1/4}\  \Gamma(\alpha n)} \   \left(x- u \right)^{\alpha n-1/2}   P_{-1/2}^{1/2-\alpha n} \left({x+u\over 2\sqrt{x u}}\right), \quad x > u > 0.\eqno(4.8)$$
Moreover, we will use the Legendre integral for the associated Legendre function \cite{erd}, Vol. I to write the right-hand side of latter equality in (4.8) in the form

$$  {\sqrt\pi  \over (xu)^{1/4}\  \Gamma(\alpha n)} \   \left(x- u \right)^{\alpha n-1/2}   P_{-1/2}^{1/2-\alpha n} \left({x+u\over 2\sqrt{x u}}\right) $$

$$=   {2   \left(x- u \right)^{2\alpha n-1}  \over \Gamma^2(\alpha n)} \   \int_0^\infty {dy \over ( 2\sqrt{xu} \ \cosh(y) + x+u)^{\alpha n}}.\eqno(4.9)$$
 Therefore, returning to (4.8) and making simple substitutions in the integral (4.9), we derive
 
 $$ { (x-u)^{2\alpha n-1} \over \Gamma(2\alpha n)} \   u^{-\alpha n}  {}_2F_1 \left( \alpha n,\ \alpha n ;\ 2\alpha n; \ 1 -{x\over u} \right) =  {2   \left(x- u \right)^{2\alpha n-1}  \over \Gamma^2(\alpha n)} \   \int_0^\infty {dy \over ( 2\sqrt{xu} \ \cosh(y) + x+u)^{\alpha n}}$$
 
 $$=  {2   \left(x- u \right)^{2\alpha n-1}  \over \Gamma^2(\alpha n)} \   \int_{(\sqrt x + \sqrt u)^2}^\infty\  {dt \over ( (t-x-u)^2-4xu )^{1/2}\  t^{\alpha n}} $$

$$=  {2   \left(x- u \right)^{2\alpha n-1}  \over \Gamma^2(\alpha n)} \  \left[  \int_{(\sqrt x + \sqrt u)^2}^{4l +1} + \int_{4l +1}^\infty \right] \  {dt \over ( (t-x-u)^2-4xu )^{1/2}\  t^{\alpha n}},\eqno(4.10) $$
where $0 < u \le x \le l$.  But since

$$  \int_{4l +1}^\infty   {dt \over ( (t-x-u)^2-4xu )^{1/2}\  t^{\alpha n}} \le   \int_{4l +1}^\infty   {dt \over ( (t-x-u)^2-4xu )^{1/2}\  t^{\alpha}} $$

$$=  (x+u)^{-\alpha} \int_{(4l +1-x-u)/ (x+u)}^\infty   {dt \over ( t^2 - (4xu)/ (x+u)^2 )^{1/2}\  (t+1)^{\alpha}} $$

$$\le   (x+u)^{-\alpha} \int_{ 2\sqrt{xu} / (x+u)}^\infty   {dt \over ( t^2 -  (4xu) / (x+u)^2 )^{1/2}\  (t+1)^{\alpha}} $$

$$=   \int_{ 1}^\infty   {dt \over ( t^2 - 1 )^{1/2}\  ( 2\sqrt{xu}\ t+x+u)^{\alpha}} \le {1\over (4 xu)^{\alpha/2} }  \int_{ 1}^\infty   {dt \over ( t^2 - 1 )^{1/2}\ t^\alpha} =  {B(1/2,\alpha/2) \over 2^{1+\alpha}  ( xu)^{\alpha/2} }, $$
and

$$  \int_{(\sqrt x + \sqrt u)^2}^{4l +1}  {dt \over ( (t-x-u)^2-4xu )^{1/2}\  t^{\alpha n}} \le x^{-\alpha n}  \int_{(\sqrt x + \sqrt u)^2}^{4l +1}  {dt \over ( (t-x-u)^2-4xu )^{1/2}} $$

$$= x^{-\alpha n}  \int_{ 2\sqrt{xu} / (x+u)}^{(4l +1 -x-u)/(x+u)}  {dt \over ( t^2- (4xu) /(x+u)^2 )^{1/2}} $$

$$=  x^{-\alpha n}\  \log\left( {4l +1 -x-u \over 2\sqrt{xu}} + \left( {(4l +1 -x-u )^2 \over 4xu} -1\right)^{1/2} \right) \le x^{-\alpha n}\  \log\left({4l+1\over \sqrt{xu}}\right),$$
we find, combining with (4.10), the following estimate  	

 $$ { (x-u)^{2\alpha n-1} \over \Gamma(2\alpha n)} \   u^{-\alpha n}  {}_2F_1 \left( \alpha n,\ \alpha n ;\ 2\alpha n; \ 1 -{x\over u} \right) $$
 
 $$ \le {2   \left(x- u \right)^{2\alpha n-1}  \over \Gamma^2(\alpha n)} \  \left[     {B(1/2,\alpha/2) \over 2^{1+\alpha} ( xu)^{\alpha/2}} +  x^{-\alpha n}\  \log\left({4l+1\over \sqrt{xu}}\right) \right].\eqno(4.11)$$
Hence, employing the H{\" o}lder and Minkowski inequalities,   the latter integral in (4.7) with  the remainder of the corresponding series can be estimated as follows

$$\left|  \int_0^x g(u) \sum_{n=N}^\infty  {\lambda^n \over \Gamma(2\alpha n)} \  (x-u)^{2\alpha n-1} u^{-\alpha n}  {}_2F_1 \left( \alpha n,\ \alpha n ;\ 2\alpha n; \ 1 -{x\over u} \right)  du\right|$$

$$\le 2 \int_0^x \left|g(u)\right|  \sum_{n=N}^\infty  {|\lambda|^n \ (x-u)^{2\alpha n-1}  \over \Gamma^2(\alpha n)}    \left[   {B(1/2,\alpha/2) \over 2^{1+\alpha}  ( xu)^{\alpha/2} } +  x^{-\alpha n}\  \log\left({4l+1\over \sqrt{xu}}\right) \right]   du$$

$$ \le 2 ||g||_{\nu,p} \left(  \int_0^x  u^{(1-\nu)p^\prime -1} \left( \sum_{n=N}^\infty {|\lambda|^n \ (x-u)^{2\alpha n-1}  \over \Gamma^2(\alpha n)}    \left[   {B(1/2,\alpha/2) \over 2^{1+\alpha}  ( xu)^{\alpha/2} } +  x^{-\alpha n}\  \log\left({4l+1\over \sqrt{xu}}\right) \right] \right)^{p^\prime}  du \right)^{1/p^\prime} $$

$$\le   {B(1/2,\alpha/2) \over 2^\alpha}\  ||g||_{\nu,p}  \sum_{n=N}^\infty {|\lambda|^n \  x^{\alpha (2n-1) -\nu}   \over \Gamma^2(\alpha n)}  \left(\int_0^1  u^{(1-\nu- \alpha/2 )p^\prime -1} (1-u)^{(2\alpha n-1)p^\prime} du\right)^{1/p^\prime} $$

$$+   2  \log\left({4l+1\over x}\right)   ||g||_{\nu,p}  \sum_{n=N}^\infty {|\lambda|^n \  x^{\alpha n  -\nu}   \over \Gamma^2(\alpha n)}  \left(\int_0^1  u^{(1-\nu )p^\prime -1} (1-u)^{(2\alpha n-1)p^\prime} du\right)^{1/p^\prime} $$

$$+  ||g||_{\nu,p}  \sum_{n=N}^\infty {|\lambda|^n \  x^{ \alpha n -\nu  }   \over \Gamma^2(\alpha n)}  \left(\int_0^1  u^{(1-\nu )p^\prime -1} (1-u)^{(2\alpha n-1)p^\prime}  \left| \log\left( u \right) \right|^{p^\prime} du\right)^{1/p^\prime} $$

$$\le    ||g||_{\nu,p}  \sum_{n=N}^\infty {|\lambda|^n \  l^{\alpha n-\nu}   \over \Gamma^2(\alpha n)} \left[ {B(1/2,\alpha/2)\  l^{\alpha n-1} \over 2^\alpha ( (1-\nu- \alpha/2 )p^\prime)^{1/p^\prime}}\  +   {2  \log\left((4l+1)/ x\right) \over ( (1-\nu )p^\prime )^{1/p^\prime}} \right.  $$

$$\left. +  {\Gamma^{1/p^\prime} (p^\prime +1) \over ((1-\nu)p^\prime)^{1+1/p^\prime} } \right] \to 0,\quad N \to \infty $$
 when $\alpha >0,\ \nu < 1-\alpha/2,\ 0 < x \le l$   owing to the rapid growth of the gamma-function.   Hence, returning to (4.7), we represent  the series inside the integral by virtue of the double hypergeometric type series  (cf.  \cite{Hai}, formula (8.2)).  Precisely, taking into account the first equality in (4.8), we  write the Gauss hypergeometric function in terms of the series to obtain

$$\sum_{n=1}^\infty  {\lambda^n \over \Gamma(2\alpha n)} \  (x-u)^{2\alpha n-1} u^{-\alpha n}  {}_2F_1 \left( \alpha n,\ \alpha n ;\ 2\alpha n; \ 1 -{x\over u} \right)$$

$$=   \lambda x^{-\alpha}  (x-u)^{2\alpha -1} \sum_{k,n =0}^\infty  \left( {\Gamma (\alpha (n+1) + k) \over \Gamma (\alpha (n+1))} \right)^2   \frac{  \Gamma(n+1) }{   \Gamma( 2\alpha (n+1)  +k)}  \left(1- {u\over x}\right)^k {\left(  \lambda x^{-\alpha}  (x-u)^{2\alpha}\right)^n\over k! \ n!}   $$

$$=   \lambda x^{-\alpha}  (x-u)^{2\alpha -1}  F \left(  \begin{matrix} ( 1, 0,1), (\alpha, 1, \alpha), (\alpha, 1, \alpha) \\  ( \alpha, 0, \alpha), (\alpha, 0, \alpha),  (2\alpha, 1, 2\alpha)\end{matrix} ;  \ 1- {u\over x}, \  \lambda x^{-\alpha}  (x-u)^{2\alpha} \right).\eqno(4.12)$$
Consequently,  the final solution (4.6) of the integral equation (4.1) has the form 

$$f(x)= g(x) +  \lambda x^{-\alpha}   \int_0^x (x-u)^{2\alpha -1}  $$

$$\times   F \left(  \begin{matrix} ( 1, 0,1), (\alpha, 1, \alpha), (\alpha, 1, \alpha) \\  ( \alpha, 0, \alpha), (\alpha, 0, \alpha),  (2\alpha, 1, 2\alpha)\end{matrix} ;  \ 1- {u\over x}, \  \lambda x^{-\alpha}  (x-u)^{2\alpha} \right)  g(u) du,\quad x \in (0,l).\eqno(4.13)$$
The results of this section are summarized in the following theorem.

{\bf Theorem 4}. {\it Let $g \in L_{\nu,p} (0,l), \ l >0,\ \nu < 1-\alpha/2,\ 1 \le p \le \infty,\ \alpha > 1/2.$ Then the Volterra-type integral equation of the second kind $(4.1),$  where $ \lambda \in \mathbb{C}:\  | \lambda| < (C_+ l^\alpha)^{-1}$  and $C_+$ is defined by $(2.10)$, has a unique solution in the space $L_{\nu,p} (0,l)$ given by formula $(4.13)$.}

\vspace{1cm}

\noindent {{\bf Acknowledgments}} \\

 \noindent The work was partially supported by CMUP [UID/MAT/00144/2019], which is funded by FCT(Portugal) with national (MEC) and European structural funds through the programs FEDER, under the partnership agreement PT2020.


\begin{thebibliography}{12}

\bibitem{But}  P.L. Butzer, A.A. Kilbas, J.J. Trujillo,  Mellin transform analysis and integration by parts for Hadamard-type fractional integrals,   
{\it J. Math. Anal. Appl.},  270 (2002), N 4,  1-15.

\bibitem{ButKil}  P.L. Butzer, A.A. Kilbas, Rodriguez- Germa, J.J. Trujillo, Stirling functions of first kind in the setting of fractional calculus and generalized differences,   
{\it Journal of Difference Equations and Applications},  13 (2007),  683-721.

\bibitem{Dattoli}  G. Dattoli, P.E. Ricci, I. Khomasuridze,  Operational methods, special polinomial and functions and solution of partial differential equations,   
{\it Integral Transforms and  Special  Functions},  15 (2004), N 4,  309- 321.

\bibitem{erd} A.  Erd{\'e}lyi, W. Magnus,  F. Oberhettinger, F. Tricomi, {\it Higher transcendental functions}. Vols. I, II. McGraw-Hill, New York (1953).

\bibitem{Hai} Nguyen Thanh Hai, S. Yakubovich,  {\it The double Mellin-Barnes type integrals and their applications to convolution theory}. Series on Soviet and East European Mathematics, 6. World Scientific Publishing,  Singapore, New Jersey  (1992). 

\bibitem{prud}  A.P. Prudnikov, Yu.A. Brychkov,  O.I. Marichev, {\it Integrals and series: Vol. 1: Elementary  Functions},
Gordon and Breach, New York  (1986); {\it Vol. 3: More Special Functions}, Gordon and Breach, New York  (1990).

\bibitem{Riordan} J. Riordan, {Combinatorial identities}, Wiley, New York (1968).


\bibitem{Samka}  S.G. Samko, A.A. Kilbas,   O.I. Marichev,  {\it Fractional integrals and derivatives. Theory and applications.} Gordon and Breach, Yverdon (1993). 

\bibitem{tit}  E.C. Titchmarsh, {\it  An Introduction to the Theory of Fourier Integrals}, Clarendon Press, Oxford ( 1937).

\bibitem{Viskov}  O.V. Viskov,  H.M. Srivastava,  New approaches to certain identities involving differential operators, 	{J. Math. Anal. Appl. } 186 (1994), 1-10. 

\bibitem{luch}  S.  Yakubovich,  Yu.  Luchko, {\it The Hypergeometric Approach to Integral Transforms and Convolutions}. Mathematics and its Applications, 287. Kluwer Academic Publishers
Group, Dordrecht (1994).





\end{thebibliography}
\end{document}